\documentclass{psapm-l}

\newtheorem{theorem}{Theorem}[section]
\newtheorem{lemma}[theorem]{Lemma}
\newtheorem{prop}[theorem]{Proposition}

\theoremstyle{definition}
\newtheorem{definition}[theorem]{Definition}
\newtheorem{example}[theorem]{Example}

\theoremstyle{remark}

\numberwithin{equation}{section}

\newcommand{\RR}{\mathbb{R}}

\newcommand{\ZZ}{\mathbb{Z}}
\newcommand{\NN}{\mathbb{N}}

\begin{document}

\title{Algebraic Recipes for Integer Programming}

\author{Bernd Sturmfels}
\address{Department of Mathematics, University of California,
Berkeley, California, 94720}
\email{bernd@math.berkeley.edu}

\thanks{This work was done while the author held
a Hewlett Packard Visiting Research Professorship 2003/2004
at MSRI Berkeley. He was also
partially supported by NSF Grant DMS-0200729}

\date{September 18, 2003.}



\begin{abstract} Integer programming is concerned with
solving linear systems of equations over the non-negative
integers. The basic question is to find a solution which
minimizes a given linear objective function for a fixed
right hand side. Here we also consider parametric versions where 
the objective function and the right hand side are allowed to vary.
The main emphasis is on Gr\"obner bases, rational generating 
functions, and how to use existing software packages.
Concrete applications to problems in statistical modeling
will be presented.
\end{abstract}

\maketitle

\section{An Introductory Coin Problem}

This lecture is about solving linear equations
over the non-negative integers.  Our point of departure is the
integer programming problem 
in \emph{standard form}:
\begin{equation}
\label{IP}
{\rm Minimize} \,\, c \cdot u \,\,\,
{\rm subject} \,\,{\rm to} \,\, A \cdot u \, = \, b \,\, 
{\rm and} \,\, u \in \NN^n.
\end{equation}
The given instance consists of an integer matrix
$A \in \ZZ^{d \times n}$, a row vector $c \in \ZZ^n$
and a column vector $b \in \ZZ^d$. The unknown
is the column vector $u = (u_1,\ldots,u_n)$.
What makes the problem  hard is the
requirement that the $u_i$ be non-negative
integers. 

As an example consider the following simple \emph{coin problem}.
Suppose you are carrying a large collection of
coins in your pocket. The allowed coins are
{\tt p}ennies ($1$ cent), {\tt n}ickels ($5$ cents),
{\tt d}imes ($10$ cents) and {\tt q}uarters ($25$ cents).
The problem is to replace your ``portfolio'' by
an equal number of coins having the same
monetary value, but such that the number of nickels 
plus the number of quarters is minimized.

This problem can be expressed in the standard form
(\ref{IP}) by setting
\begin{equation}
\label{coinIP}
d = 2,\, n = 4, \,\,\qquad
A \, = \, \begin{pmatrix}
1 & 1 & 1 & 1 \\
1 & 5 & 10 & 25 \\
\end{pmatrix} \,\,\, 
\hbox{and} \,\,\,
c = (0,1,0,1).
\end{equation}
The right hand side vector $\, b = \binom{b_1}{b_2}\,$ 
is left unspecified. Its coordinates $b_1$ and $b_2$ are the number of coins
and value (in pennies) of the portfolio respectively.

\begin{example}
\label{onefourteen}
For instance, if $\, b = \binom{10}{114}$, then we seek to express
one dollar and fourteen cents with ten coins. The optimal solution 
to this instance of (\ref{IP})
is $u = (4,2,0,4)$ and the  \emph{optimal value}
is $\, c \cdot u = 6$. In words, you can take
four pennies, two nickels and four quarters to make
one dollar and fourteen cents with ten coins, but it is impossible
to do it with less than six nickels or quarters.
\qed \end{example}

A parametric solution to our problem is provided by the 
\emph{Gr\"obner basis}
\begin{equation}
\label{gb1}
\mathcal{G} \quad = \quad \bigl\{\,\,
n^3 q - d^4, \,\, n^6 - p^5 q, \,\, n^3 d^4 - p^5 q^2, \,\,p^5 q^3 - d^8 \,\,\bigr\}. 
\end{equation}
Our Gr\"obner basis in (\ref{gb1}) is expressed as a set of \emph{monomial differences},
which is how they usually appear in computer algebra systems.
We note that there are two alternative but entirely equivalent ways 
of writing our Gr\"obner basis. In the optimization literature, it
is more common to express $\mathcal{G}$ as a  \emph{set of lattice vectors};
\begin{equation}
\label{gb2}
\mathcal{G}' \quad = \quad \bigl\{\,\,
\,(0,3,-4,1), \,\, (-5,6,0,-1),\,\, (-5,3,4,-2),\,\, (5,0,-8,3)\,\,\bigr\}.
\end{equation}

The Gr\"obner basis is a set of \emph{exchange rules} which you can
use to successively improve your  portfolio.
For instance, the first rule $\, n^3 q- d^4 \,$ says that you can replace
three nickels and one quarter with four dimes. Each of the four moves in $\mathcal{G}$
changes neither the number of coins nor their value but it decreases the
objective function value. The crucial property of being a Gr\"obner basis
says that \emph{if none of the exchange rules can be applied then
your portfolio is guaranteed to be optimal}.

There is a third way of  encoding the Gr\"obner basis,
which will be  of importance in Section 4.
Namely, we can also express $\mathcal{G}$ 
as the following \emph{generating function}:
\begin{equation}
\label{gb3}
\mathcal{G}'' \quad = \quad 
n^3 D^4 q \,+\,   P^5n^6 Q \,+\, P^5 n^3 d^4  Q^2 \,+\,p^5 D^8 q^3  .
\end{equation}
In the last representation there are two variables for each column of
the matrix $A$, and each monomial represents one exchange rule.
The lower case variable represents the gain and
the upper case variable represents the loss in the exchange of coins.

In Section 2 we explain how the Gr\"obner basis is constructed 
for an arbitrary matrix $A$ and cost function $c$,
and in Section 3 we discuss the relationship to other
notions of test sets in integer programming, including
\emph{Hilbert bases} and \emph{Graver bases}, and
we introduce Hemmecke's easy-to-use software 
\ {\tt 4ti2} \ for computing these test sets.
In Section 4 we address complexity issues. In particular,
we show how the Gr\"obner basis can computed 
in polynomial time when $d$ and $n$ are fixed.

The power of algebraic methods in integer programming stems from
the fact that they can answer parametric questions like:
\emph{What are all the optimal portfolios in our coin problem ?}
Each portfolio is given as a vector
$\, u = (u_1,u_2,u_3,u_4)\,$ or as a monomial
$\, p^{u_1} n^{u_2} d^{u_3} q^{u_4} $, and we wish
to encode all portfolios that are optimal 
solutions of (\ref{IP}) with $b = Au$. The following three
conditions are equivalent:
\begin{itemize}
\item[(a)] The vector $u \in \NN^4$ is an optimal portfolio.
\item[(b)]
None of the four monomials
$\,n^3 q $,
$\,  n^6 $,
$\, n^3 d^4 \,$ or
$\,p^5 q^3 \,$ divides $\, p^{u_1} n^{u_2} d^{u_3} q^{u_4} $.
\item[(c)]
($u_2 \! \leq \! 2 $ or $ u_4 \! = \! 0$)
and 
($u_2 \!\leq \! 5 $)
and
($u_2 \!\leq \! 2 $ or $ u_3 \! \leq \! 3$)
and 
($u_1  \! \leq \! 4 $ or $ u_4 \! \leq \! 2$).
\end{itemize}
The \emph{Hilbert series of all optimal solutions} is the
formal sum of these monomials:
$$ \sum_{u \,{\rm optimal}}
p^{u_1} n^{u_2} d^{u_3} q^{u_4} $$
This generating function is equal to the following rational function:
$$
\frac{
n^3 q^3 p^5 -n^6 q d^4 +n^6 d^4+n^3 q d^4-q^3 p^5+n^6 q-n^3 d^4-n^6-n^3 q+1}
{(1-p)(1-n)(1-d)(1-q)}.
$$
In Section 4, we will see that such Hilbert series can be computed
in polynomial time (for fixed $d$ and $n$). In Sections 5 and 6
we will focus on applications of integer programming to \emph{statistics},
and we will argue that Gr\"obner bases and generating functions
are useful tools for the optimization problems arising in this context.

The use of Gr\"obner bases as a tool for integer programming
first appeared in the paper \cite{CT} by Conti and Traverso.
Their approach was further developed in two doctoral dissertations
in the Cornell Operations Research Department, written 
by Thomas (see \cite{Tho}) and Ho\c{s}ten (see \cite{HS1})
in 1994 and 1997 respectively.
Subsequently, Ho\c{s}ten and Thomas \cite{HT1} developed an algebraic 
theory of group relaxations, extending the foundational work
in  integer programming theory which was done by 
Gomory in the 1960's. These and many other 
important topics will not be discussed in this lecture,
which aims to be introductory and self-contained. 
Readers wishing to learn more about 
commutative algebra methods in integer programming
are referred to the book \cite{Stu} and the 
survey articles \cite{HT2} and \cite{Tho2}.

\section{Gr\"obner Bases}

We consider the integer programming problem in
standard form (\ref{IP}) where $A$ and $c$ are fixed
and $b$ is arbitrary. In this section we further assume that 
$c$ is \emph{generic} in the sense that (\ref{IP}) has
a unique optimal solution for every feasible right hand side $b$.
In practise, this can always be accomplished by lexicographically
perturbing the given cost vector $c$. Consider the infinite set of all optimal solutions,
$$ {\rm Opt}_{A,c} \quad = \quad \bigl\{\,
u \in \NN^n \,: \, u \, \, \hbox{is the optimal solution of
(\ref{IP}) for }\, b = A u \, \bigr\}. $$
Suppose that $u$ and $u'$ are vectors in $\NN^n$
such that $u' \leq u$ (coordinatewise) and $u \in {\rm Opt}_{A,c}$.
Then it can be seen that $u' \in {\rm Opt}_{A,c}$. 
We paraphrase this observation in the following lemma,
using the language of partially ordered sets (= posets).

\begin{lemma} The set
$\, {\rm Opt}_{A,c}$ is an \emph{order ideal} in the partially ordered set $\,\NN^n$.
\end{lemma}

A basic result about order ideals in the poset $\NN^n$, 
known as \emph{Dickson's Lemma}, states that
the set of minimal elements in the complementary set
$\, {\rm Opt}_{A,c} \backslash \NN^n \,$ is finite.
We write  $\,{\rm Min}\bigl( \NN^n  \backslash  {\rm Opt}_{A,c}\bigr) \,$
for this finite set. Its elements are called the
\emph{minimally non-optimal points} of the integer programming family (\ref{IP}).
Recall that our introductory coin example had precisely four minimally non-optimal points:
$$ {\rm Min}\bigl( \NN^4  \backslash {\rm Opt}_{A,c}
 \bigr) \quad = \quad
\bigl\{ \,(0,3,0,1), \,\, (0,6,0,0),\,\, (0,3,4,0),\,\, (5,0,0,3)\,\,\bigr\}. $$

For every $\,g^+ \in {\rm Min}\bigl(
 \NN^n \backslash  {\rm Opt}_{A,c}\bigr) \,$
there exists a unique vector $\,g^- \in {\rm Opt}_{A,c}\,$ such
that $\, A g^+ = Ag^-$. Namely, $g^-$ is the optimal solution
to (\ref{IP}) with $\,b = A g^-$.

\begin{definition} The \emph{Gr\"obner basis} for
the matrix $A$ and cost vector $c$ is 
$$ \mathcal{G}_{A,c} \quad = \quad \bigl\{
\,\,g^+ - g^- \, \,: \,\, g^+ \in  {\rm Min}\bigl(
\NN^n \backslash {\rm Opt}_{A,c}
\bigr) \,\bigr\}. $$
This is a finite set of lattice vectors in the kernel of $A$. We
can also regard them as  monomial differences in $n$ unknowns $x_i$ or
as monomials in $2n$ unknowns $x_i, y_i$ via
$$ (2,13,0,-8,5,-7) \quad \longleftrightarrow \quad
x_1^2 x_2^{13} x_5^5 - x_4^8 x_6^7
\quad \longleftrightarrow \quad
x_1^2 x_2^{13} y_4^8 x_5^5 y_6^7  . $$
\end{definition}

The following theorem states that the Gr\"obner basis is a
minimal test set for the family of integer programs
specified by the matrix $A$ and the cost vector $c$.

\begin{theorem}
Let $u$ be a feasible solution of (\ref{IP}).
Then $u$ is non-optimal if and only if
there exists $g \in \mathcal{G}_{A,c}$ with
$g^+ \leq u$, and in this case $u-g$ is a better
feasible solution than $u$. There is no smaller set than
$ \mathcal{G}_{A,c}$ which has this property.
\end{theorem}

\begin{proof}
By construction, every element
$g = g^+ - g^-$ satisfies
$\, c \cdot g \, = \,c \cdot g^+ - c \cdot g^- > 0$.
The if-direction follows because
$\, (u-g) \cdot c < u \cdot c $ and
$g^+ \leq u$ is equivalent to $u-g$ being feasible.
For the only-if direction suppose that
no $ g \in \mathcal{G}_{A,c}$ satisfies 
$g^+ \leq u$. This means that no
element of $\,  {\rm Min}\bigl( \NN^n \backslash {\rm Opt}_{A,c}) \,$
lies below $u$ in the poset $\NN^n$. But this means
that $u \in {\rm Opt}_{A,c}$.
The minimality of $ \mathcal{G}_{A,c}$  holds
because every element of 
$\,  {\rm Min}\bigl( \NN^n \backslash {\rm Opt}_{A,c}) \,$ has to 
be reducible by some vector in the test set. \end{proof}

Under a certain genericity hypothesis on the matrix $A$,
the elements in the Gr\"obner basis are in bijection with
the \emph{neighbors of the origin}, which is
a test set for integer programming introduced by
Herbert Scarf \cite{Sca}. The connection between
neighbors and Gr\"obner bases was studied in 
a commutative algebra setting in \cite{PS}.

Let us assume that the Gr\"obner basis $ \mathcal{G}_{A,c}$ is known
to us in some explicit or implicit form. If we are given
any feasible solution  $u \in \NN^n$
then the integer programming problem (\ref{IP}) can be
solved by the following one-line algorithm:
\begin{equation}
\label{onelinealgo}
\text{
{\bf While} \ there exists $g \in  \mathcal{G}_{A,c}$ 
with $\,g^+ \leq u\,$ \ {\bf do} \
replace $u$ by $u - g$.}
\end{equation}
The problem of constructing a first
feasible solution $u$ from the right
hand side $b$ can be solved by the
same reduction process but for a different
Gr\"obner basis. The idea is completely analogous
to \emph{Phase One in the Simplex Algorithm}.
To keep our discussion simple, we will
assume that some feasible solution
$u$ is known beforehand.

One of the  objectives of this lecture is to dispel 
the belief, held by many experts in complexity theory and
combinatorial  optimization, that  the
algebraic notion of Gr\"obner bases 
is utterly useless when it comes to designing
efficient algorithms. Let me begin by pointing out that
{\bf computing Gr\"obner bases is easy and fun}.

My currently favorite tool for producing the
Gr\"obner basis $\, \mathcal{G}_{A,c}\,$ from
the matrix $A$ and the cost vector $c$ is
the software {\tt 4ti2} developed
by Raymond Hemmecke. It can be found
at {\tt www.4ti2.de} and is ridiculously easy
to download and run. It took me
(= a technologically challenged individual)
precisely three minutes  to install {\tt 4ti2}
on my (ancient) computer, and another 
minute later I was already enjoying my first
Gr\"obner basis on the screen. Actually, I don't recall having
ever encountered a piece of mathematical
software that was simpler to use than {\tt 4ti2}.

The first non-coin example I tried had $d = 3$ and $n = 7$.
The matrix $A$ was filled snakewise by 
 prime numbers and  the vector $c$ was
filled by square integers.
The input to {\tt 4ti2} consists of a matrix 
in a file named \ {\tt example} \ in the format
\begin{verbatim}
7 3
 2  3  5  7 11 13 17 
43 41 37 31 29 23 19
47 53 59 61 67 71 73 
\end{verbatim}
and a cost vector in a file \ {\tt example.cost} \
in the format
\begin{verbatim}
7 1
1 4 9 16 25 36 49
\end{verbatim}
After typing \ {\tt groebner example} \
and hitting the return key, about 
three seconds later, the  Gr\"obner basis appeared in a new
file called \ {\tt example.gro}:
\begin{verbatim}
7 241
-6 14 -11 2 1 0 0 
-6 13 -9 2 0 -1 1 
0 -1 2 0 -1 -1 1 
-4 11 -7 -1 -2 3 0 
10 -3 -35 38 -1 -7 0 
0 22 -53 39 -2 -4 0 
4 11 -46 40 0 -7 0 
 ... ... ... ...
\end{verbatim}
This Gr\"obner basis consists of $241$ vectors in $\ZZ^7$,
and it represents a solution to the parametric problem of minimizing
$\,\sum_{i=1}^7 i^2 \cdot u_i\, $ over all vectors $u \in \NN^7$ that satisfy
\begin{equation}
\label{37sys}
  \begin{pmatrix}
 2 & 3 & 5 & 7 & 11 & 13 & 17  \\
43 & 41 & 37 & 31 & 29 & 23 & 19 \\
47 & 53 & 59 & 61 & 67 & 71 & 73  \end{pmatrix} \cdot u \quad = \quad b .
\end{equation}
Knowing the $241$ vectors in the
Gr\"obner basis, we can now apply
the reduction algorithm (\ref{onelinealgo}) 
starting with any given feasible solution $u$.
Take, for instance, $\, u \, = \, (100,100,100,100,100,100)^T$.
The corresponding right hand side is 
$\,b\, =\, A \cdot u \, = \, (5800,22300,43100)$.
The algorithm in  (\ref{onelinealgo}) reduces
$u$ to the optimal solution
$\,u^* =  (62, 8, 176, 17, 423, 0, 0)$.
The optimal value is found to be
 $\,c^* =  12,525$.

Knowledge of the Gr\"obner basis allows
us to answer more advanced structural
questions about the system (\ref{37sys}).
One such question is that of finding
the \emph{integer programming gap},
a topic to be discussed in Section 4.
Another example is the question of
\emph{sensitivity analysis} with respect to the
cost function.
Suppose that the cost vector
is allowed to vary in a neighborhood
of the given vector $c$. Then the
Gr\"obner basis $\,\mathcal{G} = \mathcal{G}_{A,c}$
remains unchanged provided
$c$ ranges in the \emph{Gr\"obner cone},
which is defined by the following linear inequalities
in the unknowns $c_1,\ldots,c_n$:
\begin{equation}
\label{gbcone}
 c \cdot g > 0 \quad \hbox{for all} \,\, g 
\in \mathcal{G}. 
\end{equation}
For instance, in our coin example, the
Gr\"obner cone is the set of all solutions to
$$
3 c_2 + c_4 > 3 c_3 \,,\,\,\,
6 c_2  > 5 c_1 + c_4 \,,\,\,\,
3 c_2 + 4 c_3 > 5 c_1 + 2 c_4 \,,\,\,\,
5 c_1 + 3 c_4 > 8 c_3.
$$
The collection of all Gr\"obner cones in $\RR^n$
forms the \emph{Gr\"obner fan} of the matrix $A$.
This is an important invariant which allows us
study how the solution of (\ref{IP}) changes
as both $b$ and $c$ are allowed to vary.
See \cite{ST} for the basic theory.

The Gr\"obner fan of a matrix $A$ can be efficiently
calculated using the algorithm of Huber and Thomas
\cite{HT}. A highly efficient implementation was
recently given by Anders Jensen in his program {\tt CaTS}.
This piece of software can currently be found at the web page \
{\tt http://www.soopadoopa.dk/anders/cats/cats.html}.

\section{Hilbert bases and Graver bases}

Gr\"obner bases are closely related to
other natural notions of test sets
arising in the theory of integer programming.
A classical such notion is that of a \emph{Hilbert basis}.
Consider the problem of solving a homogeneous system
of linear equations over the non-negative integers.
As before, we assume that the defining matrix $A$ has 
$d$ rows and $n$ columns. Then our solution set is
the following semigroup:
\begin{equation}
\label{plusker}
 {\rm ker}_\NN(A) \quad = \quad \bigl\{\,
u \in \NN^n \,\, : \,\, A \cdot u = 0 \,\bigr\}. 
\end{equation}
Consider the subset of non-zero minimal elements
of the semigroup:
\begin{equation}
\label{HB}
 \mathcal{H}_A \,\,= \,\,
\bigl\{ u \in  {\rm ker}_\NN(A)\backslash \{0\} :\,
\text{no element}\,\,v \in {\rm ker}_\NN(A)\backslash \{0,u\}\,\,
\text{satisfies} \,v \leq u \bigr\}.
\end{equation}
The following result is due to
the 19th century invariant theorist Paul Gordan:

\begin{prop}
\label{gordan}
The set $\,\mathcal{H}_A\,$ is finite. It is the unique
minimal set such that
every vector in $\, {\rm ker}_\NN(A)\, $
is an $\NN$-linear combination of 
elements in $\,\mathcal{H}_A$. 
\end{prop}

The finite set $\,\mathcal{H}_A\,$ is called the \emph{Hilbert basis} of
the matrix $A$. Hilbert bases can also be computed using
the program {\tt 4ti2}. We consider the same matrix as in
(\ref{37sys}) but we now alternate the sign pattern of its columns
in the input file {\tt  example}:
\begin{verbatim}
7 3
  2  -3   5  -7  11  -13  17 
 43 -41  37 -31  29  -23  19
 47 -53  59 -61  67  -71  73 
\end{verbatim}
After typing \ {\tt hilbert example} \
and hitting the return key, about six
seconds later, the  Hilbert basis appears on a new
file called {\tt example.hil}:
\begin{verbatim}
7 1305
4 34 62 38 3 0 1 
4 35 64 38 2 1 2 
4 60 123 77 1 0 5 
4 36 66 38 1 2 3 
 ... ... ... ...
0 673 980 0 2 647 324 
0 674 982 0 1 648 325 
0 675 984 0 0 649 326 
\end{verbatim}
The Hilbert basis consists of $1,305$ vectors,
and it has a lot of internal structure.

Hilbert bases play an important role in the
recent work of Robert Weismantel and
his collaborators on ``primal methods
in integer programming''. The paper \cite{HKW}
introduces the notion of \emph{integral basis}
which is a slight generalization of Hilbert bases,
 and it presents
a simplex-like \emph{integral basis algorithm} 
which is shown to perform very well on standard benchmark
problems in integer programming.

A larger test set associated with
an integer matrix $A$ is the Graver basis,
which can be defined as follows.
For any sign pattern  $\sigma \in \{-1,+1\}^n $
let $D_\sigma$ be the $n \times n$-diagonal
matrix with $i$-th entry $\sigma_i$.
The \emph{Graver basis} of $A$ is the finite set
\begin{equation}
\label{Graverdef}
\mathcal{GR}_A \quad := \quad
\bigcup_{\sigma \in \{-1,+1\}^n }
\,D_\sigma \cdot \mathcal{H}_{A D_\sigma} 
\end{equation}
In this definition, we are taking the union over
the $2^n$ Hilbert bases for the various matrices
$A \cdot D_\sigma$. The signs are adjusted
so that each  Hilbert basis  lies in the
kernel of the original matrix $A$.
Proposition \ref{gordan}
ensures that the Graver basis $\mathcal{GR}_A$
is a finite subset of $\, {\rm ker}_\ZZ(A)$.
The following result is proved in 
\cite[\S 7]{Stu}.

\begin{prop}
\label{graverisUGB}
The Graver basis is a universal Gr\"obner basis. 
It contains, up to negating vectors, the
Gr\"obner bases of $A$ for all cost functions.
In symbols,
\begin{equation}
\label{graverisUGB2}
\bigcup_{c \in \ZZ^n}
 \mathcal{G}_{A,c}
\quad \subseteq \quad \mathcal{GR}_A .
\end{equation}
\end{prop}

The Graver basis is the ultimate test set one can 
compute for a given integer matrix $A$. It provides
a parametric solution to the integer programming
problem (\ref{IP}) when both the right hand side $b$
and the cost function $c$ are allowed to vary.

\begin{example}
\label{coingraver}
If $A$ is the $2 \times 4$-matrix (\ref{coinIP})
in our coin problem, then
$$
\mathcal{GR}_A \,\,\, = \,\,\, \bigl\{\,\,
\,(0,3,-4,1), \,\, (-5,6,0,1),\,\, (-5,3,4,-2),\,\, (5,0,-8,3)\,\,
(-5, 9, -4, 0) \bigr\}.
$$
This Graver basis has only one more element than the
Gr\"obner basis (\ref{gb2}). The advantage of the
Graver basis over the Gr\"obner basis is that we 
can now use (\ref{onelinealgo})
to solve the coin problem with respect to an
arbitrary cost vector $c$.
\qed \end{example}

The Graver basis has another natural interpretation in integer programming.
Consider our original problem (\ref{IP}) but now add the
requirement that the coordinates 
$u_i$ of the solution $u$ are bounded above
by some quantities $a_i$.
\begin{equation}
\label{IPbound}
{\rm Minimize} \,\, c \cdot u \,\,\,
{\rm subject} \,\,{\rm to} \,\, A \cdot u \, = \, b \,,
\,\, u \in \NN^n \,,\,\,
{\rm and} \,\,\,u \leq a.
\end{equation}
Here we regard $A$ and $c$ as fixed
and $(a,b) \in \ZZ^{n+d}$ as unspecified.
It turns out that the Graver basis is the unique
minimal test set for this family of integer programs.

\begin{theorem}
\label{boundtestset}
Let $u$ be a feasible solution of (\ref{IPbound}).
Then $u$ is non-optimal if and only if
there exists $g \in \mathcal{GR}_{A}$ with
$g^+ \leq u$
and $g^- \leq a-u$, and in this case $u-g$ improves
$u$. There is no smaller set than
$ \mathcal{GR}_{A}$ which has this property.
\end{theorem}

\begin{proof}
We must prove the only if direction.
Suppose $u$ is non-optimal for (\ref{IPbound})
and let $v$ be the corresponding optimal solution.
Pick $\sigma \in \{-1,+1\}^n$ so that
$\,D_\sigma (v-u)\,$ is a nonnegative vector.
There exist elements $h_1,\ldots,h_r$ in the Hilbert basis
$\,\mathcal{H}_{D_\sigma A}\,$  such that
$\, D_\sigma (v-u)\, = \, h_1 + \cdots + h_r $.
 and hence
$$ v-u \,\, = \,\, D_\sigma h_1 + \cdots + D_\sigma h_r , $$
where each summand lies in $\mathcal{GR}_A$.
Since $c \cdot (v-u) < 0$, there exists at least
one index $i$ such that $\,D_\sigma \cdot h_i \cdot c < 0$.
The vector $\, g = - D_\sigma \cdot h_i \,$ lies in 
$ \mathcal{GR}_{A}$. The construction implies that it satisfies
$\, g^+ \leq u \,$ and
$\,g^- = (D_\sigma \cdot h_i)^+ 
\leq (v-u)^+ \leq a -u$.

We now show that every element $g = g^+ - g^-$ of
$ \mathcal{GR}_{A}$ is needed in a
test set for our problem.
Suppose that $c \cdot g < 0$ and define
$a = g^+ + g^-$ and $b = A g^+ = A g^-$. 
With these choices of $a$ and $b$, the vectors
$g^+$ and $g^-$ are the only two feasible solutions 
for (\ref{IPbound}). Hence the move from $g^+$
to $g^-$ must be in the test set.
\end{proof}                                       

In light of Proposition \ref{graverisUGB} and
Theorem \ref{boundtestset}, it is highly desirable
to be able to precompute the Graver basis of a 
given integer matrix. An algorithm for this
computation is available in {\tt 4ti2}. But
the reader should be warned that 
Example \ref{coingraver} is somewhat misleading:
the Graver basis is often much larger than the
Gr\"obner basis and it takes much longer to
compute it. Consider again our 
 {\tt example} matrix,
\begin{verbatim}
7 3
 2  3  5  7 11 13 17 
43 41 37 31 29 23 19
47 53 59 61 67 71 73 
\end{verbatim}
The command \ {\tt graver example}  \
produces the Graver basis in a file
{\tt example.gra}:
\begin{verbatim}
7 29417
0 1 2 0 -1 1 1 
0 24 57 39 0 -2 2 
0 25 59 39 -1 -1 3 
 ... ... ... ...
14 -9 10 81 -6 -89 -37 
64 86 11 1 -6 -48 -28 
114 229 126 -1 -6 -11 -15 
124 268 161 -7 -6 10 -6 
 ... ... ... ...
\end{verbatim}
This Graver basis has $29,417$ elements
and it took a couple of hours to compute.
One nice feature of the Graver basis 
computation in {\tt 4ti2} is that the
program allows the exploitation of
symmetry. In many applications (e.g. in statistics)
there is a group of symmetries acting on the
columns of the matrix $A$, and the Graver basis
$\mathcal{GR}_A$ is invariant under these symmetries.
This feature allows the computation of some interesting
Graver bases whose cardinalities are in the range
of one million.

\section{The integer programming gap}

A commonly used first step towards solving a hard
integer programming problem  (\ref{IP})
is to begin by solving its \emph{linear programming relaxation}:
\begin{equation}
\label{LP}
{\rm Minimize} \,\, c \cdot u \,\,\,
{\rm subject} \,\,{\rm to} \,\, A \cdot u \, = \, b \,\,
{\rm and} \,\, u \in \RR_{\geq 0}^n.
\end{equation}
Linear programming problems are much easier
both in practise and in theory. They can be solved
in polynomial time using interior-point methods,
and the simplex algorithm performs well
in practise.
The purpose of this section is to offer algebraic
tools for comparing the hard problem 
(\ref{IP}) with the easier problem (\ref{LP}).
For an algebraic perspective on the
linear programming relaxtion see \cite{HT1}.

As before, we fix $A \in \ZZ^{d \times n}$ 
and $c \in \ZZ^n$ and regard $b \in \ZZ^d$
as unspecified. We write
$\,{\rm IPopt}_{A,c}(b)\,$ for the optimal value
of the integer program  (\ref{IP})  and we write
$\,{\rm LPopt}_{A,c}(b)\,$ for the optimal value
of the corresponding linear program  (\ref{LP}). The difference
of these quantities is a non-negative rational number
\begin{equation}
\label{gapforb}
{\rm IPopt}_{A,c}(b) \, - \,
{\rm LPopt}_{A,c}(b) \quad \geq \quad 0 .
\end{equation}
The \emph{integer programming gap} is defined as the 
maximum of the differences (\ref{gapforb})
 as $b$ ranges over
all right hand sides such that (\ref{IP}) is feasible:
\begin{equation}
\label{gap} 
{\rm gap}(A,c) \,\,\, = \,\,\,
{\rm max} \bigl\{\,
{\rm IPopt}_{A,c}(b)  - 
{\rm LPopt}_{A,c}(b) \,\,:\,\,
b \in \ZZ^d \,\,\hbox{feasible for (\ref{IP})} \,\bigr\}.
\end{equation}
It appears as if we are taking the maximum over
infinitely many different values, one for
each feasible $b$, but actually there are
only finitely many possible values 
for (\ref{gapforb}) if $A$ and $c$ are fixed, so the
maximum is attained.

\begin{example}
The integer programming gap of the coin problem
(\ref{coinIP}) equals
$$ {\rm gap}(A,c) \,\,\, = \,\,\, 76/15 \,\, = \,\, 
5.0666666... $$
This is the maximum advantage to be gained if
we allow our coins to be cut into fractional pieces.
The gap is attained for the right hand side
$\, b = \binom{10}{114} \,$ in Example \ref{onefourteen},
where $\,{\rm IPopt}_{A,c}(b) = 6\,$ is realized by
$u = (4,2,0,4)$. The optimal value of the
linear program (\ref{LP}) is 
$\,14/15 = 0.93333... \,$ and is attained by
$u' = (0,0, 136/15,14/15)$. Thus the
best way to make one dollar and fourteen
cents with ten fractional  coins is to
take $136/15$ dimes and $14/15$ quarters.
\qed \end{example}

We now give a recipe for computing the gap
by solving several auxiliary linear programming problems.
For any optimal vector $u \in \NN^n$ we define the \emph{increase set}
$$ {\rm incr}(u) \quad := \quad
\bigl\{\,i \in \{1,2,\ldots,n\} \, : \, u+e_i \,\,\hbox{optimal}\, \bigr\}. $$
A vector $u \in \NN^n$ is said to be
\emph{maximally optimal} for (\ref{IP}) if
$\,u + a $ is optimal for all vectors $a \in \NN^n$ 
whose support $\{i : a_i > 0\}$ is a subset
of the increase set $\,{\rm incr}(u)$.
For any fixed maximally optimal $u \in \NN^n$, we 
consider the
following linear program:
\begin{equation}
\label{LPmax}
{\rm Maximize} \,\, c \cdot (u-v) \,\,\,
{\rm subject} \,\,{\rm to} \,\, A \cdot (u-v) \, = \,0 \,\, 
{\rm and} \,\, v_i \geq 0 \,\,\,
\hbox{for all}\,\, i \not\in {\rm incr}(u).
\end{equation}
Here the decision variables are the
coordinates of $\,v = (v_1,\ldots,v_n)$.

\begin{theorem} \label{HoSt} {\rm (Ho\c{s}ten and Sturmfels \cite{HS2})} \
The maximum of the optimal values  of the auxiliary
linear programs (\ref{LPmax}), as $u$ ranges over all
maximally optimal solutions to (\ref{IP}),
coincides with the integer programming gap, $\,{\rm gap}(A,c)$.
\end{theorem}

\begin{example}
Our coin problem (\ref{coinIP}) has three maximally optimal solutions:
\begin{equation}
\label{threemaxopt}
 \bigl(4,2,\underline{0},\underline{0} \bigr) \, , \,\,\,
\bigl(\underline{0},2,\underline{0},2 \bigr) \,,\,\,
\bigl(\underline{0},5,3,0 \bigr).  
\end{equation}
In each case the increase set ${\rm incr}(u)$ is indicated
by the underlined coordinates.
The vectors in (\ref{threemaxopt}) are easily derived from the
Gr\"obner basis (\ref{gb1}).
For instance, the last portfolio
(consisting of five nickels and three dimes)
is maximally optimal because
adding one nickel, dime or quarter makes
that portfolio non-optimal but
adding any number of pennies
is fine. The program (\ref{LPmax})
for that portfolio equals
\begin{eqnarray*}
& {\rm Maximize} \,\,\,
 4 - v_2 - v_4  \,\,\,
 \hbox{subject to} \,\,\,
\, \, v_1 + v_2 + v_3 + v_4 = 4 \\
& \quad
 v_1 + 5 v_2 + 10 v_3 + 25 v_4 = 60 \, , \,\,\,
v_2 \geq 0 \,\,\,\, \hbox{and} \,\,\,\,  v_4 \geq 0. 
\end{eqnarray*}
The optimal value is $4$. The
optimal values of (\ref{LPmax}) for
$\, (4,2,\underline{0},\underline{0} ) \, $ and
$\, (\underline{0},2,\underline{0},2) \,$
are $76/15$  and $5$ respectively,
and hence  $\,
{\rm gap}(A,c) \, = \,  {\rm max}\bigl\{
4,5,76/15\} \, = \, 76/15$.
\qed \end{example}

Theorem \ref{HoSt}
 furnishes an algorithm for computing
$\, {\rm gap}(A,c)\,$ because the set of
maximally optimal solutions to (\ref{IP})
is always finite and can be computed from the
Gr\"obner basis $\, \mathcal{G}_{A,c}\,$
by the algebraic process of
\emph{irreducible decomposition of monomial ideals}.
A highly efficient implementation of this process
was developed  by Alex Milowski  in his
Master's thesis project at
San Francisco State University.  The non-trivial
gap computations in Examples \ref{37hard} and \ref{K5gap}
were done by Ho\c{s}ten and Milowski
 using {\tt 4ti2} (to derive the
Gr\"obner basis),  Milowski's software 
(to get the maximally optimal solutions)
and {\tt maple} (to solve the linear programs (\ref{LPmax})).

\begin{example}
\label{37hard}
Let $d = 3$, $n = 7$ and  consider the instance
discussed  in (\ref{37sys}):
$$ A \, = \,
  \begin{pmatrix}
 2 & 3 & 5 & 7 & 11 & 13 & 17  \\
43 & 41 & 37 & 31 & 29 & 23 & 19 \\
47 & 53 & 59 & 61 & 67 & 71 & 73  \end{pmatrix},\,\,\quad
c = (1,4,9,16,25,36,49). $$
There are $553$ maximally optimal solutions, and the gap is
\begin{equation}
\label{234gap}
 {\rm gap}(A,c) \quad = \quad
43771/183 \,\,\, = \,\,\, 239.1857923
\end{equation}
The gap is attained by the
right hand side
$\, b \, = \,(661, 1710, 3994)^T$.
For this choice of $b$, the optimal value of
(\ref{IP}) equals $1,757$ and is given by
the optimal solution $\, u = (7, 4, 0, 22, 0, 3, 26)$,
while the optimal value of
(\ref{LP}) is a little less than $1,518$ and is given
by the optimal solution 
$\,u =  (0,0,0,0, 14029/244, 463/366, 521/732)$.
\qed
\end{example}

\section{Short rational generating functions}

The importance of rational generating functions
for lattice  point problems has been known to
combinatorialists  for a long time.  Their role 
as an efficient tool in integer programming, however,
has been recognized only quite recently, in response
to the polynomial time  algorithms of
Barvinok \cite{Bar} and Barvinok-Woods \cite{BW}.
This work was further extended by
De Loera et.al.~\cite{DHHHSY}, \cite{DHHHY}.
This section reports on these methods and their
implementations in the software {\tt LattE}.

As a point of entry consider the following  variant of our problem:
\emph{List all optimal solutions to the integer program} (\ref{IP}).
For a concrete example take $d=1,\,n=4, \, 
A = (\, 1 \,\,\, 1 \,\,\, 1 \,\,\, 1 \,) $, $c= (0,0,0,1)$,
and suppose $\, b \gg 0$. Here (\ref{IP}) equals
\begin{equation}
\label{mini4}
 \hbox{Minimize} \,\, u_4\,\,
\hbox{subject to} \,\,\, 
u_1 + u_2 + u_3 + u_4 \, = \, b \,\,
\hbox{and} \,\,\,
u_1 , u_2 , u_3 , u_4 \in \NN. 
\end{equation}
The set of optimal solutions is the set of all lattice points 
$(u_1,u_2,u_3,0)$ in a large triangle. We can 
write them all down as the terms of the generating function
\begin{equation}
\label{gf4}
\sum_{u \,{\rm optimal} \,{\rm for} \,(\ref{mini4})} 
\!\!\!\!
x_1^{u_1}
x_2^{u_2}
x_3^{u_3}
x_4^{u_4} \quad = \quad
\sum_{u_1=0}^b 
\sum_{u_2=0}^{b-u_1}
x_1^{u_1} x_2^{u_2} x_3^{b-u_1-u_2}.
\end{equation}
The number of terms in this series equals
$\,(b+1)(b+2)/2 = O(b^2)$. This quantity is
exponential in the size of the input, which
is $O({\rm log}(b))$. Indeed, the number of
bits needed to write down the  line
(\ref{mini4}) grows like the logarithm
of the integer $b$, while the number
of terms on the right hand side of (\ref{gf4}) is exponential in
${\rm log}(b)$.
It appears to be impossible to ``list''
all feasible solutions to  (\ref{mini4})
in polynomial time, given that their
number grows exponentially in the input size.
Nonetheless, it can be done, namely, by
rewriting (\ref{gf4}) as the  \emph{short rational 
generating function}
$$
x_1^{b} \cdot
(1-\frac{x_2}{x_1})^{-1} 
(1-\frac{x_3}{x_1})^{-1} 
\,\, + \,\,
x_2^{b} \cdot
(1-\frac{x_1}{x_2})^{-1} 
(1-\frac{x_3}{x_2})^{-1} 
\,\, + \,\,
x_3^{b} \cdot
(1-\frac{x_1}{x_3})^{-1} 
(1-\frac{x_2}{x_3})^{-1} .
$$
The reader is invited to check that this rational function
equals the series (\ref{gf4}).
The rational function can be computed 
in time $\, O({\rm log}(b))\,$
and it represents the ``list'' of all optimal
solutions to (\ref{mini4}). This approach works
for any integer program:

\begin{theorem} 
Suppose that $d$ and $n$ are fixed. Then
the number of  optimal solutions to {\rm  (\ref{IP})}
and the rational generating function 
$\, \sum \{\, x^u \, : \, u \,\, {\rm optimal}\,\,
{\rm for} \,\, (\ref{IP})\, \} $,
which encodes the  set of optimal solutions,
can be computed in polynomial time. 
\end{theorem}

\begin{proof}
The optimal value $c^*$ of (\ref{IP})
can be computed in polynomial time
using Lenstra's algorithm \cite{Len}. 
Now apply Barvinok's 
lattice point algorithm \cite{Bar}
to the polytope
$\, \{ \, u \in \RR_{\geq 0}^n \,:\,
A u = b, \, c \cdot u = c^* \,\}$. It computes
the desired
generating function
and its evaluation at $(1,1,\ldots,1)$ in polynomial time.
\end{proof}

The techniques underlying Barvinok's algorithm
were developed substantially further by
Barvinok and Woods \cite{BW}. Using their 
\text{Projection Theorem}, one can derive
polynomial-time algorithms based on
rational generating functions for
essentially all of the algorithmic
questions we have encountered so far.
We refer to \cite{BW} and
\cite{DHHHSY} for proofs
of various parts of
the following theorem.

\begin{theorem}
\label{polytime}
Consider a matrix $A \in \ZZ^{d \times n}$ and
a vector $c \in \ZZ^d$ whose dimensions
$d$ and $n$ are fixed. Then the rational
generating functions which encode
the following sets can be computed
in time  polynomial in the bit complexity of $A$ and $c$:
\begin{enumerate}
\item the Gr\"obner basis $\, \mathcal{G}_{A,c}$,
\item the set $\, {\rm Opt}_{A,c} \,$ of all optimal solutions,
\item the set
$\,{\rm Min}\bigl( \NN^n  \backslash  {\rm Opt}_{A,c}\bigr) \,$
of minimally non-optimal points,
\item the Hilbert basis $\, \mathcal{H}_A$,
\item the Graver basis $\, \mathcal{GR}_A$,
\item the set of maximally optimal solutions, and
\item the integer programming gap $\, {\rm gap}(A,c)$.
\end{enumerate}
\end{theorem}

The result (7) about the gap appears in \cite{HS2}.
The objects in (1)-(6) are highly structured
subsets of $\ZZ^n$. It is this special structure
which allows for a short encoding. For encoding
the Gr\"obner basis $\, \mathcal{G}_{A,c}$,
the paper \cite{DHHHSY} uses
a generating function in $2n$ variables
as in  (\ref{gb3}). But all the sets
in (1)-(6) can also be coded as
formal sums of Laurent monomials
(representing vectors in $\ZZ^n$),
and the Barvinok-Woods method will
give short rational functions 
for these encodings.

\vskip .1cm

A magnificent computer program for solving
lattice point problems by means of short rational generating
functions has been developed
by the group of Jesus De Loera at UC Davis.
It is called {\tt LattE} and can be obtained at the web site
{\tt http://www.math.ucdavis.edu/$\sim$latte/}. This program
can be used to count the number of feasible solutions
to an integer program (\ref{IP}) as follows.

Consider our coin problem in (\ref{coinIP})
with $\,b = \binom{999}{5000}$, so we wish to 
arrange  $999$ coins to be worth fifty dollars.
In order to determine in how many ways this
can be accomplished, we create the following
{\tt LattE} input file which we call {\tt coins}:
\begin{verbatim}
6 5
999 -1  -1  -1  -1 
5000 -1 -5 -10 -25 
 0    1  0   0   0
 0    0  1   0   0
 0    0  0   1   0
 0    0  0   0   1
2 1 2
\end{verbatim}
The command \ {\tt latte equ coins} \
will count the number of feasible solutions 
to (\ref{IP}).
The output which appears on the screen reveals that
the answer is $9,352$:
\begin{verbatim}
This is LattE v1.0 beta.   (September 17, 2002)
Revised version.           (Aug        1, 2003)
The polytope has 4 vertices.
 ....  ....  .... ... 
Creating generating function.
Starting final computation.

  **** THE GRAND TOTAL IS: 9352 ****

Computation done. 
Time: 0.01 sec
\end{verbatim}
This run of {\tt LattE} has created
the following output on a new file called
{\tt coins.maple}:
\begin{verbatim}
gF:=x[0]^4999170*x[1]^(-8000506)*x[2]^1000*x[3]^166/((1-x[0]^(-4995)
*x[1]^7993*x[2]^(-1)*x[3]^(-1))*(1-x[0]^5*x[1]^(-9)*x[3]^(-1)))+
x[0]^4999170*x[1]^(-8000506)*x[2]^1000*x[3]^166/((1-x[0]^(-9985)
  ....   ....    ....   ....   ....   ....    ....   .... 
... *x[1]^7993*x[2]^(-1)*x[3]^(-1))*(1-x[0]^(-5)*x[1]^9*x[3]));
\end{verbatim}
This is the short rational generating function representing
the formal sum of $ 9,352$ monomials, one for each feasible
solution. You can get the expanded form of this
generating function
reading the file 
\ {\tt coins.maple} \ into the computer
algebra system {\tt maple}. After you have done this,
please type the {\tt maple} command
\ {\tt simplify(gF);}.

The program {\tt LattE} can also be used
to solve the minimization problem
(\ref{IP}). To this end
we need to add the cost vector
\ {\tt 0 1 0 1} \, in 
an extra line at the end of the
input file \ {\tt coins}.
Typing now the {\tt LattE} command sequence
{\tt ./latte equ min coins2}, we obtain the following 
output on the screen
\begin{verbatim}
 ....  ...    ....    ....
All cones have been decomposed.
6 cones in total.
Computing the points in the Parallelepiped of the unimodular Cones.
An optimal solution for [0 1 0 1] is: [555 2 441 1].
The optimal value is: 3.
The gap is: 3
Computation done.
\end{verbatim}
We conclude from this {\tt LattE} session that the best way of making 
fifty dollars with $999$ coins is to take
one quarter, two nickels, $441$ dimes and $555$ quarters.

I  tried {\tt LattE} on considerably bigger problems
and I found that it performs quite well. The speed
is particularly impressive for knapsack problems
(d = 1) with large integer coefficients. 
For this class of problems, {\tt LattE} is faster
than the current version of the commerical software
{\tt CPLEX} on some
instances. This parallels the observation,
already made in \cite{HS1}, that programs like {\tt CPLEX}
are not always the best choice
for low-dimensional problems with large integers,
given that they are designed for highly structured
0/1 problems with many variables.

The authors of {\tt LattE}  informed me that they
intend to incorporate
all of the tasks listed in Theorem \ref{polytime}
into a future version of their program.
The lesson to be learned here is that
algebraic software like {\tt 4ti2} and {\tt LattE}
can definitely play a useful role in the box of tools
available to practitioners
of integer programming.

\section{Some integer programs arising in statistics}

We present an application to the
statistical theory of {\em disclosure limitation}.
See  \cite{CG} and \cite{DF} and the references  therein.
Suppose we are given data in the form of an
$n$-dimensional table of nonnegative integers. The aim is
to release some marginals of the table but not the
table's entries themselves. If the range of possible values that a
particular entry can attain in any table satisfying the released
marginals is too narrow then this entry may be
exposed. This shows the importance of determining tight  upper
and lower bounds for each entry in a given table. 

A choice of marginals corresponds to  fixing subsets
$\,F_1, \ldots, F_k\,$ of $\{1, \ldots, n\}$. It  can be represented by
a  zero-one matrix $A$, as described in
 \cite[\S 1]{HSu}. In statistical language, the matrix
$A\,$ specifies a {\em hierarchical model} for a 
{\em contingency table} with $n$ {\em factors}.
Suppose $v$ is a table with nonnegative integer entries,  where
  the marginals are computed according to a fixed hierarchical model $A$
  and let $v_{i_1 i_2 \cdots i_n}$
be a particular cell of the table $v$.
What we are interested in is the following
 \textbf{table entry security problem}: \emph{Compute
optimal lower and upper bounds $L$ and $U$ such that $L
  \leq u_{i_1 i_2 \cdots i_n} \leq U$ for all tables 
$u$ which have the same  marginals as $v$}.

The table entry security problem is an integer program:
minimize (or maximize)  $u_{i_1 i_2 \cdots i_n}$
 over all tables with nonnegative {\em
  integer} entries subject to fixing the marginals.
In order to write this integer program
in the standard form (\ref{IP}), we need to give
the precise definition of the relevant matrices $A$.
Consider $d_1 \times \cdots \times d_n$-tables with 
entries $u_{i_1 i_2 \cdots i_n}$ where $1 \leq i_j \leq d_j$.
We fix a  hierarchical model by specifying  $\,F_1, \ldots, F_k $.
The marginals of our table
are computed with respect to these subsets.
If $F_i = \{j_1, \ldots, j_s\}$ then the
$F_i$-marginal is a $d_{j_1} \times \cdots \times d_{j_s}$ table 
$b$ with  entries 
\begin{equation}
\label{umumum}
   b_{k_1 \cdots  k_s} \quad = 
    \sum_{i_{j_1} = k_1, \ldots, i_{j_s} = k_s} u_{i_1 \cdots i_n}. 
\end{equation}

We define $A$ to be the zero-one matrix with
$d_1 d_2 \cdots d_n$ columns
representing the linear map that computes 
the marginals of tables. 
We let $u$ be  the vector of variables representing the cell entries.
Then $A \cdot u$ represents the $k$ lower-dimensional
tables computed as in (\ref{umumum}).  
The \textbf{table entry security problem} is
\begin{equation}
\label{newip}
{\rm Minimize} \,\, ({\rm Maximize})
  \,\,\,\, u_{11 \cdots  1} \quad
\hbox{subject to} \quad A \cdot u  =  b, 
\, u \geq 0, \, u \,\, \hbox{integral}.
\end{equation}

Here we only consider the cell entry $u_{1 1 \cdots 1} $
(corresponding to the first column of $A$) because
there is a
transitive symmetry group acting on the columns of $A$.

\begin{example}
The classical {\em transportation problem} \cite[p.~221]{Sch}
corresponds
to $d_1 \times d_2$-tables where the marginals are computed
with respect to $F_1 = \{1\}$ and $F_2 = \{2\}$. The 
{\em three-dimensional transportation problem} \cite{Vla}
concerns $d_1 \times d_2 \times d_3$-tables with $F_1 = \{1,2\}$, 
$F_2 = \{1,3\}$, and $F_3 = \{2,3\}$. The marginals are
$$ b_{ij} = \sum_k u_{ijk}\, ,\, \,\,\, 
 b_{ik} = \sum_j u_{ijk} \,, \,\,\,\,
b_{jk} = \sum_i u_{ijk}. $$      
For a discussion from the Gr\"obner basis perspective see
\cite[\S 14.C]{Stu}. \qed
\end{example}

\begin{example} Consider the \emph{four-cycle model}
for binary random variables. Here $n = 4$,
$d_1 = d_2 = d_3 = d_4 = 2$,
$F_1 = \{1,2\}$, 
$F_2 = \{2,3\}$, 
$F_3 = \{3,4\}$, and  $F_4 = \{1,4\}$.
The matrix $A\,$ has 
$\, d_1d_2 + d_2d_3 + d_3d_4 + d_1d_4 \, = \,16\,$ rows
and it has
$\,d_1d_2d_3d_4 = 16$ columns. 
We write it in {\tt 4ti2} format on a file
name {\tt fourcycle}:
\begin{verbatim}
16
16
 1 1 1 1 0 0 0 0 0 0 0 0 0 0 0 0
 0 0 0 0 1 1 1 1 0 0 0 0 0 0 0 0 
 0 0 0 0 0 0 0 0 1 1 1 1 0 0 0 0
 0 0 0 0 0 0 0 0 0 0 0 0 1 1 1 1 
 1 1 0 0 0 0 0 0 1 1 0 0 0 0 0 0
 0 0 1 1 0 0 0 0 0 0 1 1 0 0 0 0  
 0 0 0 0 1 1 0 0 0 0 0 0 1 1 0 0  
 0 0 0 0 0 0 1 1 0 0 0 0 0 0 1 1  
 1 0 0 0 1 0 0 0 1 0 0 0 1 0 0 0  
 0 1 0 0 0 1 0 0 0 1 0 0 0 1 0 0  
 0 0 1 0 0 0 1 0 0 0 1 0 0 0 1 0  
 0 0 0 1 0 0 0 1 0 0 0 1 0 0 0 1  
 1 0 1 0 1 0 1 0 0 0 0 0 0 0 0 0 
 0 1 0 1 0 1 0 1 0 0 0 0 0 0 0 0  
 0 0 0 0 0 0 0 0 1 0 1 0 1 0 1 0
 0 0 0 0 0 0 0 0 0 1 0 1 0 1 0 1
\end{verbatim}
The Graver basis of this matrix $A$
consists of $106$ vectors. The
{\tt 4ti2} command \ {\tt graver fourcycle} \
delivers the Graver basis on a new file
{\tt fourcycle.gra}:
\begin{verbatim}
16 106
0 0 0 0 0 1 0 -1 1 -1 0 0 -1 0 0 1 
-1 0 1 0 1 0 -1 0 1 0 0 -1 -1 0 0 1 
0 0 0 0 0 1 0 -1 0 0 0 0 0 -1 0 1 
0 0 0 0 1 0 -1 0 0 0 0 0 -1 0 1 0 
1 0 -1 0 0 0 0 0 -1 0 1 0 0 0 0 0 
0 0 -1 1 0 0 1 -1 0 0 0 0 0 0 0 0 
-1 1 0 0 1 -1 0 0 0 0 0 0 0 0 0 0 
0 -1 1 0 0 1 -1 0 1 0 0 -1 -1 0 0 1 
  ... ... ... ...  ... ... ... ...
\end{verbatim}
In light of Proposition \ref{graverisUGB},
we can use this Graver basis with 
the Algorithm (\ref{onelinealgo})
to solve (\ref{IP}) for any
cost function. In particular, we can use it
solve (\ref{newip}).
\qed
\end{example}

As the parameters $d_1,\ldots,d_n$ increase,
it becomes harder to solve the
integer program (\ref{newip}) exactly.
Researchers in disclosure limitation
have resorted to solving the
linear programming relaxation (\ref{LP}) instead:
 minimize (or maximize)  $u_{i_1 i_2 \cdots i_n}$ over all
tables with nonnegative {\em real} entries subject to fixing the
marginals. This relaxation is tractable,  but it usually fails
 to deliver the exact integers
$L$ and $U$.  One faces the problem of finding the 
integer programming gap for the table entry security problem. 
This application was the original motivation
for the paper \cite{HS2}.

\begin{example} 
\label{K5gap}
What follows may serve as a test case
for future software for computing the gap.
We consider the $K_5$-model for five binary
random variables. Here
$\,n = 5, \,k = 10 ,\,
d_1 = \cdots = d_5 = 2$ and the
$F_i$ are the ten two-element subsets of $\{1,2,3,4,5\}$.
The cost function is $\, c = - e_{11111}$, 
corresponding to maximizing in (\ref{newip}).
The matrix $A$ has $40$ rows and $32$ columns,
and it has rank $ 16 $. We found 
\begin{equation}
\label{thegapisthree}
 {\rm gap}(A,c) \quad = \quad 3. 
\end{equation}
The gap is attained by the following $2 \times 2 \times 
2 \times 2 \times 2$-table:
\begin{equation}
\label{opttable}
 u \quad = \quad
e_{11112} + e_{11121} +
2 \cdot e_{11211} + 2 \cdot e_{12111} + 2 \cdot e_{21111} 
+ e_{22222}. 
\end{equation}
This table is optimal for the maximization
problem in  (\ref{newip}). The $40$ entries of the
 right hand side vector $b$ are the entries in the
 ten marginal $2 \times 2$-tables:
$$
\begin{pmatrix}
u_{11***} & u_{12***} \\
u_{21***} & u_{22***} \\
\end{pmatrix}
\,\, = \,\,
\begin{pmatrix}
 4 & 2  \\
 2 & 1 \\
\end{pmatrix}
\, ,\,\, \ldots \, \ldots \, , \,\,
\begin{pmatrix}
u_{***11} & u_{***12} \\
u_{***21} & u_{***22} \\
\end{pmatrix}
\,\, = \,\,
\begin{pmatrix}
 5 & 1  \\
 2 & 1 \\
\end{pmatrix}
$$
Since the unit table $e_{11111}$ does not 
appear in (\ref{opttable}), we
have $\,{\rm IPopt}_{A,c}(b)\,= \, 0 $.

The optimal value of the linear programming
relaxation equals $\,{\rm LPopt}_{A,c}(b)\,= \, 3 $.
This value is attained by the following
 fractional $2 \times 2\times 2 \times 2 \times 2$-table
\begin{eqnarray*}
v \quad = &
 3 \cdot e_{00000} +
\frac{2}{3} e_{11211} + 
\frac{1}{3} e_{11221} + 
\frac{2}{3} e_{12111} + 
\frac{1}{3} e_{12122} + 
\frac{1}{3} e_{12211} + 
\frac{1}{3} e_{12221} + \\ &
\frac{1}{3} e_{12222} + 
\frac{2}{3} e_{21112} + 
\frac{1}{3} e_{21121} + 
\frac{2}{3} e_{21211} + 
\frac{1}{3} e_{21222} + 
\frac{2}{3} e_{22111} + 
\frac{1}{3} e_{22121} .
\end{eqnarray*}
We invite the reader to check
that the tables $u$ and $v$
have the same marginals.
\qed
\end{example}

\bibliographystyle{amsalpha}

\end{document}